\documentclass[12pt]{article}

\usepackage{amsmath}

\usepackage{amsfonts}

\usepackage{pstricks}

\numberwithin{equation}{section}

\title{Perfect Quadrilateral Right Prisms}
\author{Allan J. MacLeod,\\Statistics and Mathematics Group,\\
University of the West of Scotland,\\High St., Paisley,\\Scotland.  PA1 2BE\\
(e-mail: allan.macleod@uws.ac.uk)}
\date{}

\begin{document}

\maketitle

\begin{abstract}
We consider right prisms with horizontal quadrilateral bases and tops, and vertical rectangular sides. We look for examples where
all the edges, face diagonals and space diagonals are integers. We find examples when the base is an isosceles trapezium and a parallelogram,
but no solution for a kite or rhombus.
\end{abstract}

\section{Introduction}
Let {\bf ABCD} be a convex quadrilateral in the x-y plane and let {\bf EFGH} be the identical shape
vertically translated $h$ units. Then the $3$-dimensional object with vertices {\bf ABCDEFGH} is known as a right prism.
In this paper, we consider finding such prisms where the lengths of the edges, face diagonals and space diagonals
can all be integers. If {\bf ABCD} and {\bf EFGH} are rectangles. then the prism is in fact a cuboid, and the problem
is that of finding a perfect rational cuboid, which is still a classic unsolved problem, see Bremner \cite{br1}, Dickson \cite{dick1},
Guy \cite{guy}, Leech\cite{leech1} and van Luijk \cite{luijk}.

Very recently, Sawyer and Reiter \cite{sawreit} have shown the existence of perfect parallelipipeds. These are not normally right prisms, in that
the solutions found do not have rectangular vertical faces. In this paper, we consider quadrilaterals such as the
trapezium, the parallelogram, the kite and the rhombus, to see whether perfect right prisms with these base shapes exist. It is interesting
to note that a mathematician of the calibre of Kummer \cite{kumm} considered the problem of quadrilaterals with integer sides and diagonals.

\section{Basics}
Since the sides of a right prism are clearly rectangles, there will be several equations of the form
\begin{displaymath}
a^2 + b^2 = \Box \, , \hspace{1cm} a, b \in \mathbb{Z}
\end{displaymath}
which have a general solution $a = \alpha (m^2-n^2), b=\alpha 2mn$. The $\alpha$ term can be a nuisance so we get rid
of it by reforming the solution as $a/b=P(r)=(r^2-1)/(2r)$ with $r=m/n$.

Often, an equation in a problem reduces to the equation
\begin{equation}
P(c)=d
\end{equation}
where we wish rational solutions. A simple consideration of the underlying quadratic shows that
\begin{equation}
1+d^2=\Box
\end{equation}
is a necessary condition.

The most usual such occurrence is when we look for rational solutions of
\begin{displaymath}
x^2+h^2=\Box \hspace{1cm} , \hspace{1cm} y^2+h^2=\Box
\end{displaymath}

Thus, $x/h=P(a)$ and $y/h=P(b)$ for some rational $a$ and $b$. Then, $P(b)=(y/x)P(a)$, which has a rational solution
if $1+y^2P^2(a)/x^2=\Box$, which implies
\begin{equation}
d^2=y^2a^4+(4x^2-2y^2)a^2+y^2
\end{equation}
has rational solutions. There has clearly at least one rational solution, when $a=0, d=\pm y$, so the quartic is birationally equivalent to an elliptic curve,
see Silverman and Tate \cite{siltate}.

Using standard transformations given in Mordell \cite{mord}, we find the equivalent curve to be
\begin{equation}
v^2=u^3+(x^2+y^2)u^2+x^2y^2u=u(u+x^2)(u+y^2)
\end{equation}
with $a=v/(y(u+x^2))$.

The curve has torsion subgroup (usually) isomorphic to $\mathbb{Z}_2 \times \mathbb{Z}_4$, with finite torsion points of order $2$ at
$(0,0)$, $(-x^2,0)$, $(-y^2,0)$, and points of order $4$ at $(xy \, , \, \pm xy(x+y))$, and $(-xy \, , \, \pm xy(x-y))$. None of these points give a non-trivial value of
$a>1$. We thus need curves which have rank greater than $0$. For example, $x=2,y=3$ leads to $v^2=u^3+13u^2+36u$ which has rank $0$, and
hence there are no non-trivial solutions of $4+h^2=\Box, \; 9+h^2=\Box$.

For $x=5,y=2$, however, $v^2=u^3+29u^2+100u$ has rank $1$ with generator $(2,18)$.
This point gives $a=1/3$ which is not acceptable. If we add $(2,18)$ to one of the order $4$ torsion points, we get the point $(-20,40)$, which
gives $a=4$, giving $P(a)=15/8$ and $h=8/3$. Scaling $x,y,h$ by $3$ gives $x=15,y=6,h=8$.

We look for as many generators as we can easily find. Let $G_1,\dots,G_s$ be the independent generators found with $0 \le s \le r$
where $r$ is the rank of the curve, then we compute the points
\begin{equation}
Q = n_1 G_1 + n_2 G_2 + \ldots + n_s G_s + T_k
\end{equation}
where $-L \le n_i \le +L$ and $T_k$ is a torsion point, $0 \le k \le 7$, where $T_0$ is the point at infinity. For each point $Q$ we determine $a, P(a)$ and
$h=x/P(a)$.

Suppose $P$ is a point of infinite order. We find $P$ and $P+T$ where $T$ is a point of order $2$ generates the $4$ values of $\pm a_1,\pm 1/a_1$, for some $a_1$,
as we use the $3$ points of order $2$ and the point at infinity. For the four values from $P+T$ where $T$ varies over the $4$ points of order $4$, we get $\pm a_2, \pm 1/a_2$.
It is clear, however, that these values all give the same value for $|P(a_1)|$ and $|P(a_2)|$, and hence we only generate $2$ different positive values of $h$. Further,
we have the product of these $h$ values equals $|x\,y|$. All this can be proven by standard application of the formulae for addition on elliptic curves, using a
good symbolic algebra package to alleviate the boredom.

To find whether a curve has strictly positive rank, for a given $(x,y)$ pair, we use the Buhler-Gross-Zagier \cite{bgz} method of computing the L-series
and its derivatives to estimate the rank, if the conductor of the curve is not too large. If the rank is estimated at least $2$ or rank=1 and the Birch \& Swinnerton-Dyer
conjecture predicts a small height for the generator, we search for integer points on the minimal model of the curve, which are then transformed to possible
generators of the original curve.

\section{General Quadrilateral}
This is a convex quadrilateral with no parallel sides and no equal or complementary angles, where we can assume w.l.o.g. that $x$
is the largest side.

 \begin{center}
\setlength{\unitlength}{1cm}
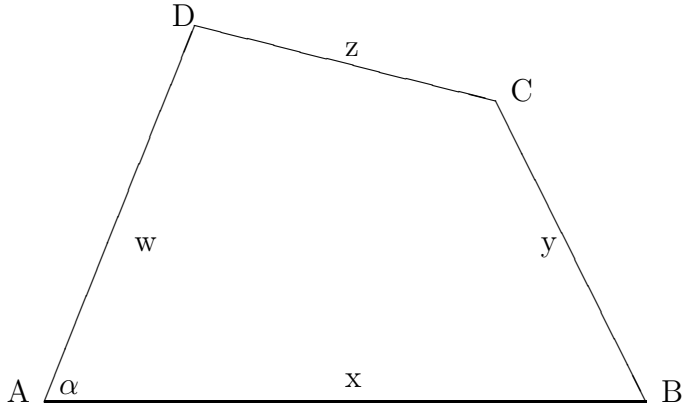
\begin{picture}(9,5)
\put(0,0){\line(2,5){2}}
\put(2,5){\line(4,-1){4}}
\put(6,4){\line(1,-2){2}}
\put(0,0){\line(1,0){8}}
\put(-0.5,0){A}
\put(1.7,5){D}
\put(8.2,0){B}
\put(6.2,4){C}
\put(4,0.2){x}
\put(4,4.6){z}
\put(1.2,2){w}
\put(6.6,2){y}
\put(0.2,0.1){$\alpha$}
\end{picture}
\end{center}

We generate values for $x$ and $w$ and $\cos \alpha=m/n$, since integer diagonals imply rational cosines.  We first test if $BD=v$ is a rational solution of
$v^2=x^2+w^2-2xw\cos \alpha$. If it is, we generate sides $y,z$ in the range $v< y+z < 2x$ and determine $\cos \angle ABD$,
$\sin \angle ABD$, $\cos \angle DBC$ and $\sin \angle DBC$, by elementary trigonometry.

\begin{center}
TABLE 1\\Near-solution lengths for general quadrilaterals\\
\begin{tabular}{rrrrrrrrr}
$\;$&$\;$&$\;$&$\;$&$\;$ \\
x&y&z&w&u&v&h&$\cos \alpha$&$\cos B$\\
  91 & 80 & 45 & 63 & 45 & 35 & 60 & 25/26 & 113/130 \\
 128 & 56 & 36 & 56 & 88 & 88 & 105 & 23/28 & 23/28 \\
 171 & 150 & 84 & 150 & 192 & 192 & 80 & 29/100 & 29/100 \\
 159 & 147 & 48 & 147 & 171 & 171 & 140 & 37/98 & 37/98 \\
 273 & 112 & 385 & 273 & 343 & 273 & 180 & 1/2 & -1/2 \\
441 & 210 & 96 & 210 & 294 & 294 & 280 & 23/28 & 23/28 \\
2223 & 2185 & 513 & 2185 & 2432 & 2432 & 420 & 9/23 & 9/23 \\
7735 & 7616 & 3927 & 7616 & 9401 & 9401 & 3960 & 1/4 & 1/4 \\
23205 & 16779 & 5712 & 16779 & 20349 & 20349 & 2860 & 49/94 & 49/94
\end{tabular}
\end{center}

These relations give $\cos B$, and we test the rationality of $u=AC$ from $u^2=x^2+y^2-2xy\cos B$. If we find rational $u$ and $v$, we test for
\begin{displaymath}
\{x, y, z, w, u, v\}^2 + h^2 = \Box
\end{displaymath}
by using the elliptic curve method of section $2$. We find, if possible, values of $h$ with $\{x^2,y^2\}+h^2=\Box$ and test these with
$\{z,w,u,v\}$.

Despite extensive testing we were unable to find a perfect prism but we found several where only one of the side or space diagonals
is not an integer. These are given in Table 1. A quick survey of Table $1$ shows that all the quadrilaterals have at least one pair of lengths identical. We have been unable
to find any quadrilaterals, with all horizontal lengths distinct, where more than $4$ vertical diagonals are rational. The $5$ results found so far are given
in Table $2$.

\begin{center}
TABLE 2\\Distinct-sided quadrilaterals with $4$ rational vertical diagonals.\\
\begin{tabular}{rrrrrrrrr}
$\;$&$\;$&$\;$&$\;$&$\;$ \\
x&y&z&w&u&v&h&$\cos \alpha$&$\cos B$\\
 152&140&56&36&88&154&105&1/16&23/28\\
 390&380&315&95&374&425&168&-5/19&251/475\\
 420&60&600&288&472&540&175&-2/15&-191/225\\
 728&504&462&266&640&735&2520&5/32&25/49\\
 1672&616&396&1540&1664&968&1155&23/28&167/847
\end{tabular}
\end{center}

\section{Isosceles Trapezium}
An isosceles trapezium is a shape as shown in the diagram, where $AB$ and $DC$ are parallel, and
$x,y,z \in \mathbb{Z}$. Thus $\cos\alpha=(x-y)/2z$, so that $\cos \alpha$ is rational.

The perpendicular gap $w$ satisfies
\begin{equation}
w^2=z^2-(x-y)^2/4
\end{equation}
and the base diagonal is $v=AC=BD$, which satisfies
\begin{equation}
v^2=\frac{(x+y)^2}{4}+w^2=z^2+xy
\end{equation}

\begin{center}
\setlength{\unitlength}{1cm}
\begin{picture}(9,5)
\put(0,0){\line(1,2){2}}
\put(2,4){\line(1,0){4}}
\put(6,4){\line(1,-2){2}}
\put(0,0){\line(1,0){8}}
\put(-0.5,0){A}
\put(1.7,4){D}
\put(8.2,0){B}
\put(6.2,4){C}
\put(4,0.2){x}
\put(4,3.6){y}
\put(1.2,2){z}
\put(6.6,2){z}
\put(0.2,0.1){$\alpha$}
\put(7.6,0.1){$\alpha$}
\put(2,0){\line(0,1){4}}
\put(2.2,2){w}
\end{picture}
\end{center}

The vertical side diagonals give the three equations
\begin{equation}
x^2+h^2=\Box \hspace{0.5cm} , \hspace{0.5cm} y^2+h^2=\Box \hspace{0.5cm} , \hspace{0.5cm} z^2+h^2=\Box
\end{equation}
whilst there is only one space diagonal equation, namely $v^2+h^2=\Box$.

We have $x/h=P(a)$, $y/h=P(b)$, $z/h=P(c)$
for some rational $a,b,c >1$. We, thus, generate values of $a,b,c$, find $P(a), P(b), P(c)$ and hence
suitable integer $x, y, z, h$. We reject those for which, firstly $z^2+xy \ne \Box$ and then $v^2+h^2 \ne \Box$.

This simple search was coded in Pari-GP and quickly found several solutions. The $10$ smallest (using the measure $\|\cdot\|=x+y+z+v+h$)
are given in Table $3$.

\begin{center}
TABLE 3\\Solution lengths for isosceles trapezia\\
\begin{tabular}{ccccc}
$\;$&$\;$&$\;$&$\;$&$\;$ \\
x&y&z&v&h\\
  364& 275& 320& 450& 240\\
 1152& 507& 780& 1092& 1040\\
 3325& 1053& 1620& 2475& 2160\\
 3328& 361& 1824& 2128& 3420\\
 3549& 2601& 1326& 3315& 1768\\
 4225& 2527& 2223& 3952& 2964\\
 5632& 4693& 1368& 5320& 1824\\
 2754& 1984& 4455& 5031& 7560\\
 6647& 3168& 5950& 7514& 2040\\
 10633& 2300& 4550& 6720& 5040
\end{tabular}
\end{center}

We can parameterise the lengths of the isosceles trapezium by parameterising $v^2=x^2+z^2-2xz\cos\alpha$. Standard
algebra leads to the fomulae,
\begin{equation}
x=2(c+d) \hspace{0.4cm} , \hspace{0.4cm} y=2d(1+c\,d) \hspace{0.4cm} , \hspace{0.4cm} z=1-d^2 \hspace{0.4cm} , \hspace{0.4cm} v=2cd+1+d^2
\end{equation}
where $c,d < 1$, and, in fact, $c=\cos\alpha \equiv i/j$ with $i,j \in \mathbb{Z}$.

The expression for $z$ leads to the temptation to define $h=2d$ so that $z^2+h^2=(1+d^2)^2$. Given this $h$, for $x^2+h^2$ and $y^2+h^2$ to be squares, we
need, respectively
\begin{equation}
c^2+2cd+2d^2=\Box \hspace{1cm} , \hspace{1cm} c^2d^2+2cd+2=\Box
\end{equation}
while $v^2+h^2$ square needs
\begin{equation}
d^4+4cd^3+(4c^2+6)d^2+4cd+1=\Box
\end{equation}

The quadric $c^2+2cd+2d^2=\Box$ can be easily parameterised by $c=g^2-2f^2$ and $d=2f(f-g)$, but when we substitute these into the second quadric we get
a sextic in $g$ or an octic in $f$ which must be made square. These curves are likely to have genus greater than $1$ and so only a finite
number of rational solutions, by Mordell's theorem.

The latter quartic clearly is birationally equivalent to an elliptic curve and we find this curve to be
\begin{equation}
v^2=u^3+i^2u^2-j^4u \hspace{1cm} , \hspace{1cm} d=\frac{v-iu}{j(u+j^2)}
\end{equation}

This curve seems, from numerical experiments just to have the single finite torsion point $(0,0)$. It has the rational points $(j^2,\pm ij^2)$ and
$(-j^2,\pm ij^2)$. If we double the point $(j^2,ij^2)$ we get $(j^4/i^2,-j^6/i^3)$ and since we assume $i/j$ is a completely
reduced fraction, we have non-integer coordinates, so by the Nagell-Lutz theorem this is of infinite order.

Using $u=j^4/i^2, v=j^6/i^3$ gives
$d=(j^3-i^2j)/(i^3+ij^2)$ which means that $x^2+h^2=\Box$ when $i^8+3i^4j^4-2i^2j^6+2j^8=\Box$, whilst $y^2+h^2$ is square when $i^4+2i^2j^2+5j^4=\Box$.
The latter is birationally equivalent to the elliptic curve $v^2=u^3-u^2-u$ which has rank $0$, so this choice of point does not lead to all
vertical diagonals square.

%

%

\section{Cyclic Quadrilateral}
The isosceles trapezium is an example of a cyclic quadrilateral - a quadrilateral whose $4$ vertices lie on
the circumference of a circle.

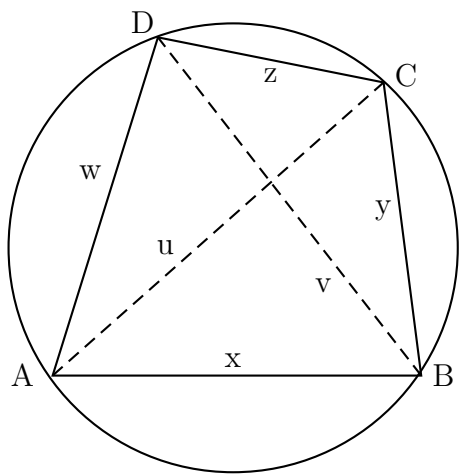
\begin{center}
\setlength{\unitlength}{1cm}
\begin{pspicture}(9,9)
\pscircle(4,4){3}
\pspolygon(1.6,2.3)(6.5,2.3)(6,6.2)(3,6.8)
\psline[linestyle=dashed](1.6,2.3)(6,6.2)
\psline[linestyle=dashed](6.5,2.3)(3,6.8)

\rput(1.2,2.3){A}
\rput(6.8,2.3){B}
\rput(6.3,6.3){C}
\rput(2.8,7.0){D}

\rput(4,2.5){x}
\rput(6,4.5){y}
\rput(4.5,6.3){z}
\rput(2.1,5){w}
\rput(3.1,4){u}
\rput(5.2,3.5){v}

\end{pspicture}
\end{center}

We look for lengths $(x,y,z,w,u,v)$ which are all positive integers first. This has been considered by Schubert \cite{schub} and
Sastry \cite{sast1} \cite{sast2}.
There is an error in the basic formulae in \cite{sast1}, whilst \cite{sast2} uses a somewhat unusual parameterisation.

We derive a more standard form here. The results are based on the well-known formula for the length of an arc of a circle of radius $R$, namely
$2 R \sin \theta$ where $\theta$ is the angle the arc makes with the circumference. Let $\alpha=\angle ADB=\angle ACB$, $\beta=\angle BAC=\angle BDC$,
$\gamma=\angle CAD=\angle CBD$ and $\delta=\angle ABD=\angle ACD$. Since ABCD is a cyclic quadrilateral, opposite angles are complimentary,
thus $\alpha+\beta+\gamma+\delta=\pi$, so there are only $3$ free choices of angle.

Therefore, assuming w.l.o.g. than $2R=1$, we have
\begin{equation}
x=\sin \alpha \hspace{1cm} , \hspace{1cm} y=\sin \beta \hspace{1cm} , \hspace{1cm} z =\sin \gamma \hspace{1cm}
\end{equation}
\begin{equation}
w=\sin\alpha\cos\beta\cos\gamma+\cos\alpha\sin\beta\cos\gamma+\cos\alpha\cos\beta\sin\gamma-\sin\alpha\sin\beta\sin\gamma
\end{equation}
\begin{equation}
u=\sin\alpha\cos\beta+\cos\alpha\sin\beta \, , \, v=\sin\beta\cos\gamma+\cos\beta\sin\gamma
\end{equation}

Thus, the sines and cosines of these angles are rational so we can assume
\begin{equation}
\cos \alpha=\frac{f^2-1}{1+f^2} \, \, , \, \, \sin\alpha=\frac{2f}{1+f^2} \hspace{0.2cm} , \hspace{0.2cm}
\cos \beta=\frac{g^2-1}{1+g^2}  \, \, , \, \,\sin\beta=\frac{2g}{1+g^2} \hspace{0.2cm} , \hspace{0.2cm}
\end{equation}
\begin{displaymath}
\cos \gamma=\frac{h^2-1}{1+h^2} \, \, , \, \, \sin\gamma=\frac{2h}{1+h^2} \hspace{0.2cm} , \hspace{0.2cm}
\end{displaymath}
with $f,g,h \in \mathbb{Q}$.

Clearing denominators, we have
\begin{equation}
x=f(g^2+1)(h^2+1) \hspace{1cm} , \hspace{1cm} y=g(f^2+1)(h^2+1)
\end{equation}
\begin{equation}
z=h(f^2+1)(g^2+1) \hspace{0.5cm} , \hspace{0.5cm} w=(f(g+h)+gh-1)(f(gh-1)-g-h)
\end{equation}
\begin{equation}
u=(f+g)(fg-1)(1+h^2) \hspace{1cm} , \hspace{1cm} v =(g+h)(gh-1)(1+f^2)
\end{equation}

Despite extensive computing using these formulae, we have been unable to generate a cyclic quadrilateral, with all sides
and diagonals mutually distinct, where more than $3$ of the $t^2+h^2$ terms are square. The smallest quadrilateral with $3$
squared terms found so far is $x=561$, $y=750$, $z=1275$, $w=1560$, $u=1575$, $v=1197$, $h=400$, $\cos \alpha=11/2$,
$\cos \beta =4$ and $\cos \gamma=2$.

\section{Parallelogram}
The parallelogram $ABCD$ has sides $AB \parallel DC$ and $AD \parallel BC$, with $AB=DC=x$ and
$AD=BC=y$.
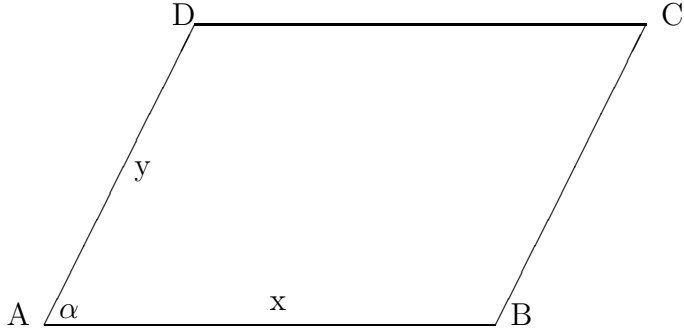
\begin{center}
\setlength{\unitlength}{1cm}
\begin{picture}(9,5)
\put(0,0){\line(1,2){2}}
\put(2,4){\line(1,0){6}}
\put(6,0){\line(1,2){2}}
\put(0,0){\line(1,0){6}}
\put(-0.5,0){A}
\put(1.7,4){D}
\put(6.2,0){B}
\put(8.2,4){C}
\put(3,0.2){x}
\put(1.2,2){y}
\put(0.2,0.1){$\alpha$}
\end{picture}
\end{center}

Let $AC=z$ and $BD=w$. Then
\begin{equation}
z^2=(x+y \cos \alpha)^2+(y \sin \alpha)^2=x^2+y^2+2xy\cos \alpha
\end{equation}
\begin{equation}
w^2=(x-y \cos \alpha)^2+(y \sin \alpha)^2=x^2+y^2-2xy \cos \alpha
\end{equation}

These are both homogeneous quadrics, and the intersection of two such elements
is known to be an elliptic curve. We can parameterise the first equation as
\begin{displaymath}
\frac{x}{y}=\frac{2(\cos \alpha - \mu)}{\mu^2-1}
\end{displaymath}
which we substitute into the second quadric giving
\begin{equation}
d^2=\mu^4+4k\mu^3+2(1-2k^2)\mu^2-12k\mu+8k^2+1
\end{equation}
where $k=\cos \alpha=m/n$ since $\cos \alpha$ clearly rational.

This quartic is birationally equivalent to an elliptic curve, which we eventually find to be
\begin{equation}
h^2=g^3+2(m^2+n^2)g^2+(m^2-n^2)^2g
\end{equation}
with the relation
\begin{equation}
\mu=\frac{h-2m(m^2-n^2)}{n(g-(m^2-n^2))}
\end{equation}

The elliptic curve has torsion subgroup isomorphic to $\mathbb{Z}_2 \times \mathbb{Z}_4$, with finite torsion points
of order $2$ for $g=0,-(m+n)^2,-(m-n)^2$ and order $4$ points at $g=\pm (m^2-n^2)$. None of these lead to non-trivial $\mu$ and
$x/y$. Thus we need curves of strictly positive rank. We use the same approach as described in section $2$ to find curves
with infinite generators. We generate multiples of such generators and these give different $\mu$ and hence different
$x,y,z,w$.

For given $x,y$ values, we use the method of section $2$ to find values of $h$, if any, where $x^2+h^2=\Box$ and $y^2+h^2=\Box$. These values of
$h$ are then tested for $z^2+h^2=\Box$ and $w^2+h^2=\Box$.

Unfortunately, despite extensive testing, no solutions were found. In fact, only rarely did we have one of $z^2+h^2$ or $w^2+h^2$ a square,
and this almost always occurred when $z$ or $w$ equalled $x$ or $y$.

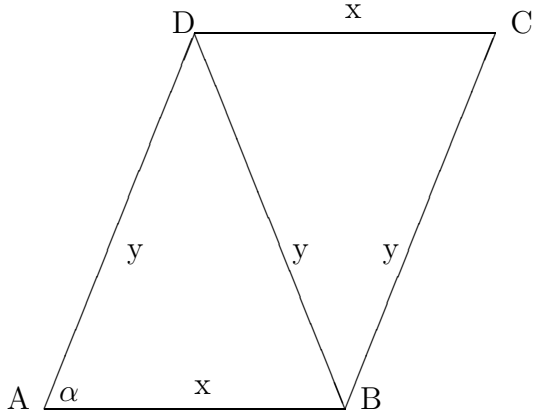
\begin{center}
\setlength{\unitlength}{1cm}
\begin{picture}(9,6)
\put(0,0){\line(2,5){2}}
\put(2,5){\line(1,0){4}}
\put(4,0){\line(2,5){2}}
\put(0,0){\line(1,0){4}}
\put(2,5){\line(2,-5){2}}
\put(-0.5,0){A}
\put(1.7,5){D}
\put(4.2,0){B}
\put(6.2,5){C}
\put(2,0.2){x}
\put(4,5.2){x}
\put(1.1,2){y}
\put(4.5,2){y}
\put(3.3,2){y}
\put(0.2,0.1){$\alpha$}
\end{picture}
\end{center}

Consider the special parallelogram $ABCD$ given above. Then $\cos \alpha=x/(2y)$ and $AC^2=z^2=2x^2+y^2$. This can be simply parameterized
by $x=2pq, y=p^2-2q^2, z=p^2+2q^2$. We can thus easily generate values of $x,y,z$ and hence use section $2$ to find values of $h$.

A large amount of testing has found a single solution where all sides, face diagonals and space diagonals are integers. This occurs where $x=2522, y=16199$,
giving $w=16199$, $z=16587$, and $\cos \alpha = 13/167$. If we have a height $h=9240$, we have
side diagonals of $9578$ and $18649$, and the two space diagonals $18649$ and $18987$.

It might be thought that we can set $h=p^2-q^2$ since $x^2+h^2=\Box$. For $y^2+h^2$ and $z^2+h^2$ to be square we need, respectively,
$2p^4-6p^2q^2+5q^4=\Box$ and $2p^4+2p^2q^2+5q^4=\Box$. Both have obvious solutions, so are birationally equivalent to elliptic curves.
We find the curves and transformations to be
\begin{equation}
E_1: v^2=u^3+12u^2-4u \hspace{0.5cm} , \hspace{0.5cm} \frac{p}{q}=\frac{4u+v+4}{2u+v-4}
\end{equation}
\begin{equation}
E_2: v^2=u^3-4u^2-36u \hspace{0.5cm} , \hspace{0.5cm} \frac{p}{q}=\frac{4u+v+12}{v-2u-12}
\end{equation}
in the correct order.

Both curves have only one finite torsion point $(0,0)$, both have rank $1$ with $E_1$ having generator $G_1=(1,3)$ and $E_2$ generator $G_2=(9,9)$.
We can easily generate, using Pari-GP, the points on $E_1$, compute $p/q$ and see if the second quartic is square. We have done this
up to $(u,v)=99G_1$ with no luck.

\section{Kite}
A standard kite has the shape shown where the diagonals meet at right-angles with $BD$ bisected by $AC$. Call $\angle DAC=\alpha$ and $\angle DCA=\beta$.

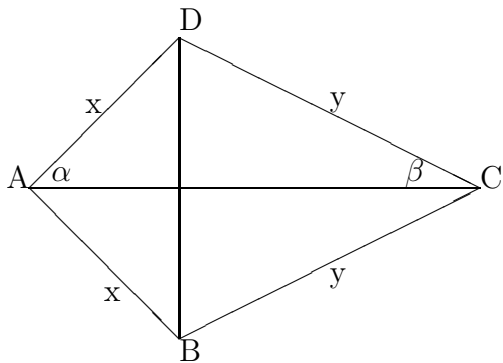
\begin{center}
\setlength{\unitlength}{1cm}
\begin{picture}(8,6)
\put(1,3){\line(1,1){2}}
\put(1,3){\line(1,-1){2}}
\put(3,1){\line(2,1){4}}
\put(3,5){\line(2,-1){4}}
\put(1,3){\line(1,0){6}}
\put(3,1){\line(0,1){4}}
\put(0.7,3){A}
\put(3,0.7){B}
\put(3,5.1){D}
\put(7,3){C}
\put(2,1.5){x}
\put(1.75,4){x}
\put(5,1.75){y}
\put(5,4.1){y}
\put(1.3,3.1){$\alpha$}
\put(6.0,3.1){$\beta$}
\end{picture}
\end{center}

Thus $AC=z=x\cos \alpha +y \cos \beta$ and $BD=w=2x\sin \alpha=2y\sin \beta$. Thus we must have the sines and cosines
rational, so we assume
\begin{displaymath}
\cos \alpha=\frac{f^2-1}{f^2+1} \; , \; \sin \alpha =\frac{2f}{kf^2+1} \; , \; \cos \beta=\frac{g^2-1}{g^2+1} \; , \; \sin \beta =\frac{2g}{g^2+1},
\end{displaymath}
where $f,g \in \mathbb{Q}$ and $f,g > 1$.

The equation for $w$ gives
\begin{equation}
g^2-2kg+1=0 \hspace {1cm} , \hspace{1cm} k=\frac{y}{x \sin\alpha}
\end{equation}
so that a rational $g$ needs $k^2-1=\Box$, and $g=k+\sqrt{k^2-1}$.

We first find $x$ and $y$ such that the elliptic curve from section $2$ has positive rank. We then generated values of $f=m/n$ with $m+n\le299$, and
tested for rational $g$. If such $g$ existed, we computed $z$, $w$ and $h$ and tested for $z^2+h^2=\Box$ and $w^2+h^2=\Box$. In the testing so far, no
example with both conditions satisfied has been found. Many examples with one satisfaction were found, the vast majority with $z=x$ or $z=y$.
Table $4$ lists some examples where all the lengths are different.

\begin{center}
TABLE 4\\Nearly-perfect kites\\
\begin{tabular}{rrrrrrr}
$\;$&$\;$&$\;$&$\;$&$\;$&$\;$&$\;$ \\
x&y&z&w&h&$\cos\alpha$&$\cos\beta$\\
   75& 435& 450& 144& 308& 7/25& 143/145\\
  100& 208& 252& 160& 105& 3/5& 12/13\\
  585& 1190& 1375& 1008& 1200& 33/65& 77/85\\
   3300& 13260& 13800& 6336& 9625& 7/25& 1073/1105
\end{tabular}
\end{center}

It should be noted that it is easy to see that we cannot have $w=x$ or $w=y$.

We can simply parameterise a kite as having sides and diagonals which are multiples of
\begin{equation}
x=g(f^2+1) \hspace{0.2cm} , \hspace{0.2cm} y=f(g^2+1) \hspace{0.2cm} , \hspace{0.2cm} z= f(g^2-1)+g(f^2-1) \hspace{0.2cm} , \hspace{0.2cm} w=4fg
\end{equation}

If we define $h=4f^2-g^2$ then $w^2+h^2=\Box$. To have $x^2+h^2=\Box$, we must have
\begin{equation}
(g^2+16)f^4-6g^2f^2+g^2(g^2+1)=\Box
\end{equation}
and a little searching shows that $f=g/2$ gives a solution. Thus this quartic is birationally equivalent to an elliptic curve.

Standard analysis shows that the equivalent curve is
\begin{equation}
v^2=u^3+12p^2q^2u^2-4p^2q^2(p^2+4q^2)^2u
\end{equation}
where $g=p/q$ with $p,q$ integers with no common factors. Numerical investigations suggest that the obvious point $(0,0)$ is the only common torsion point,
and that, in general, the torsion subgroup is isomorphic to $\mathbb{Z}_2$.

The transformation back from the elliptic curve is
\begin{equation}
f=\frac{p(v+4q^2u+8p^2q^2(p^2+4q^2))}{2q(v-p^2u-8p^2q^2(p^2+4q^2))}
\end{equation}

If we set the denominator to $0$, we get $v=p^2(u+8p^2q^2+32q^4)$, and if we substitute into the elliptic curve, we find
$u=-16p^2q^2, v=\pm 8p^2q^2(p^2-4q^2)$. Since we suspect that there is only one finite torsion point, this would
have to be a point of infinite order. Doubling this point gives
\begin{displaymath}
u=\frac{(p^4+72p^2q^2+16q^4)^2}{16(p^2-4q^2)^2}
\end{displaymath}
and the fractional nature of this shows that $-16p^2q^2$ will give a point of infinite order for most $(p,q)$ pairs. Thus, the
elliptic curve has rank at least $1$ in most cases, with the largest rank found so far being $4$, starting with $p=17, q=8$.
 So, we can easily find kites with $(x^2,w^2)+h^2=\Box$.

We wrote a simple program which searched for other generators on these curves, generated points on the curve, and
hence $(x,y,z,w,h)$ sets which we checked for more than two $a^2+h^2=\Box$. These were very sparsely distributed and we only
found one new solution with $3$ squares - $(x,y,z,w,h)$ = $(8745, 4182, 10881, 6336, 14840)$.

It might be thought that making $y^2+h^2=\Box$ first would produce a similar outcome. We require,
\begin{equation}
16f^4+(g^4-6g^2+1)f^2+g^4=\Box
\end{equation}
which can be linked to the elliptic curve
\begin{equation}
v^2=u^3-2(p^4-6p^2q^2+q^4)u^2+(p^2+q^2)^2(p^2+4pq+q^2)(p^2-4pq+q^2)u
\end{equation}
with $g=p/q$ and transformation $f=v/(8q^2u)$.

These curves can have rank $0$ and have torsion subgroup isomorphic to $\mathbb{Z}_2 \times \mathbb{Z}_4$. Using this form,
we found no perfect kites and only one new near-perfect kite with sides $(1883,1924,1107,3640,2400)$.

Finally, we can look to make $(w^2,z^2)+h^2=\Box$. Using the prior parameterisations, we require
\begin{equation}
(g^2+16)f^4+2g(g^2-1)f^3+(g^4-12g^2+1)f^2-2g(g^2-1)f+g^2(g^2+1)=\Box
\end{equation}
where, again, we find $f=g/2$ gives a solution. Transforming, we find the elliptic curve
\begin{equation}
v^2=u^3+(p^4+18p^2q^2+q^4)u^2+12p^2q^2(p^4-4q^4)u
\end{equation}
where $g=p/q$ and the inverse transformation is
\begin{equation}
f=\frac{p(v+(p^2-5q^2)u+24p^2q^2(p^2-2q^2))}{2q(v+(q^2-2p^2)u-24p^2q^2(p^2-2q^2))}
\end{equation}

As before, the torsion subgroup seems to only contain $(0,0)$ as a finite point.
Solving either the numerator or denominator for $v$ and substituting into the elliptic curve gives the
point $(-16p^2q^2,\pm 8p^2q^2(p^2+4q^2))$, so the curves have rank at least $1$. Numerical experiments show that these curves
can have quite a large rank - the largest found is $6$ starting with $p=29, q=15$ - and the rank was always at least $2$.
This led to investigating the simple integer points and the discovery of a second point of infinite order
$(-(p^2+q^2)^2,\pm2pq(p^2+q^2)(p^2+4q^2))$. Surprisingly, given all points available, we were unable to find any kites with even $3$ vertical diagonals integer.

\section{Rhombus}
A rhombus is both a kite and a parallelogram with shape
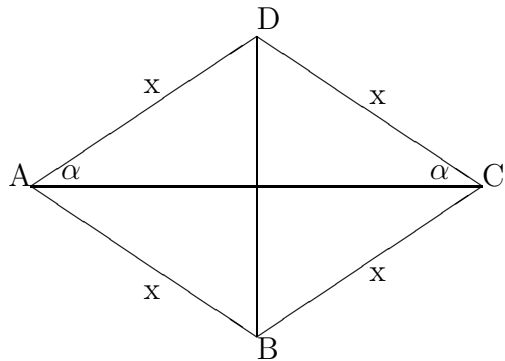
\begin{center}
\setlength{\unitlength}{1cm}
\begin{picture}(8,6)
\put(1,3){\line(3,2){3}}
\put(1,3){\line(3,-2){3}}
\put(4,1){\line(3,2){3}}
\put(4,5){\line(3,-2){3}}
\put(1,3){\line(1,0){6}}
\put(4,1){\line(0,1){4}}
\put(0.7,3){A}
\put(4,0.7){B}
\put(4,5.1){D}
\put(7,3){C}
\put(2.5,1.5){x}
\put(2.5,4.25){x}
\put(5.5,1.75){x}
\put(5.5,4.1){x}
\put(1.4,3.1){$\alpha$}
\put(6.3,3.1){$\alpha$}
\end{picture}
\end{center}

Let $AC=z=2x\cos\alpha$ and $BD=w=2x\sin\alpha$, so $\cos\alpha$ and $\sin\alpha$ must be rational so set
\begin{displaymath}
\cos\alpha=\frac{f^2-1}{f^2+1} \hspace{1cm} , \hspace{1cm} \sin\alpha=\frac{2f}{f^2+1}
\end{displaymath}
where $f=m/n>1$.

We thus need
\begin{equation}
x^2+h^2=\Box \hspace{0.25cm} , \hspace{0.25cm} \frac{4(m^2-n^2)^2x^2}{(m^2+n^2)^2}+h^2=\Box \hspace{0.25cm} , \hspace{0.25cm}
\frac{16m^2n^2x^2}{(m^2+n^2)^2}+h^2=\Box
\end{equation}
and we can determine the simple parameterisation
\begin{equation}
x=m^2+n^2 \hspace{0.6cm} , \hspace{0.6cm} z=2(m^2-n^2) \hspace{0.6cm} , \hspace{0.6cm}
w=4mn
\end{equation}

If we set $x=p^2-q^2$ and $h=2pq$, we have $x^2+h^2=\Box$. Thus we need $m^2+n^2+q^2=p^2$. One parameterisation is
$m=2ab, \, n=2ac, \, q=a^2-b^2-c^2, \, p=a^2+b^2+c^2$. Substituting, we find $z^2+h^2=\Box$ requires
\begin{equation}
4(a^8+2(7b^4-18b^2c^2+7c^4)a^4+(b^2+c^2)^2)=\Box
\end{equation}
whilst $w^2+h^2=\Box$ needs
\begin{equation}
4(a^8-2(b^4-30b^2c^2+c^4)a^4+(b^2+c^2)^2)=\Box
\end{equation}

It might be thought, at first, that this is an octic expression and therefore fairly unmanageable. Note, however,
that $x=4a^2(b^2+c^2), \, z=8a^2(b^2-c^2), \, w=16a^2bc$, so the term of interest is actually $a^2$, which
makes these expressions quartics and hence equivalent to elliptic curves.

We find that the expression for $z^2+h^2$ is related to the curve
\begin{equation}
v^2=w^3-(7b^4-18b^2c^2+7c^4)w^2+4(b^2-c^2)^2(b^2-3c^2)(3b^2-c^2)w
\end{equation}
and $w^2+h^2$ related to
\begin{equation}
v^2=w^2+(b^4-30b^2c^2+c^4)w^2-16b^2c^2(b^4-14b^2c^2+c^4)w
\end{equation}
with, in both cases, $a^2=v/w$. Both curves have at least $8$ torsion points, with, for most $(b,c)$ pairs,
the torsion subgroup isomorphic to $\mathbb{Z}_2 \times \mathbb{Z}_4$.

Extensive testing using both elliptic curves to generate values of $a^2$ from input $(b,c)$ pairs has
not found any examples of a perfect right rhombus prism. For $x^2,z^2+h^2=\Box$ a simple example is
$x=75, z=42, w=144, h=40$, and , for $x^2,w^2+h^2=\Box$ a simple example is $x=988, z=760, w=1824, h=315$.

\end{document}